\documentclass[11pt]{amsart}  
\usepackage{amsmath,amssymb,amsfonts,hyperref,times}

\title{Lower bounds for measurable chromatic numbers}

\author
{Christine Bachoc} 
\address{C.~Bachoc, Laboratoire A2X, Universit\'e Bordeaux I, 351,
cours de la Li\-b\'e\-ration, 33405 Talence, France}
\email{bachoc@math.u-bordeaux1.fr}

\author {Gabriele Nebe} \address{G.~Nebe, Lehrstuhl D f\"ur
  Mathematik, RWTH Aachen University, Templergraben 64, 52062 Aachen,
  Germany} \email{nebe@math.rwth-aachen.de}

\author
{Fernando M\'ario de Oliveira Filho}
\address{F.M.~de Oliveira Filho, Centrum voor Wiskunde en Informatica (CWI),
Kruislaan 413, 1098 SJ Amsterdam, The Netherlands}
\email{f.m.de.oliveira.filho@cwi.nl}

\author
{Frank Vallentin} 
\address{F.~Vallentin, Delft Institute of Applied Mathematics, Technical University of Delft, P.O. Box 5031, 2600 GA Delft, The Netherlands}
\email{f.vallentin@tudelft.nl}

\thanks{The third author was partially supported by CAPES/Brazil under
grant BEX 2421/04-6. The fourth author was partially supported by the Deutsche
Forschungsgemeinschaft (DFG) under grant SCHU 1503/4.}


\date{July 17, 2009}

\newcommand{\defi}[1]{\textit{#1}}

\newcommand{\Ccal}{\mathcal{C}}
\newcommand{\R}{\mathbb{R}}

\newcommand{\chim}{\chi_{\textup{m}}}

\newtheorem{defin}{Definition}[section]

\newtheorem{proposition}[defin]{Proposition}
\newtheorem{theorem}[defin]{Theorem}
\newtheorem{remark}[defin]{Remark}
\newtheorem{corollary}[defin]{Corollary}

\newtheorem{example}[defin]{Example}

\DeclareMathOperator{\ort}{O}
\DeclareMathOperator{\Aut}{Aut}
\DeclareMathOperator{\Harm}{Harm}

\newcommand{\intsn}{\int_{S^{n-1}}\int_{S^{n-1}}}

\newcommand{\al}{\alpha}
\newcommand{\be}{\beta}

\newcommand{\paak}{P_k^{(\al,\al)}}
\newcommand{\paaj}{P_j^{(\al,\al)}}
\newcommand{\paajp}{P_{j+1}^{(\al,\al)}}
\newcommand{\pabk}{P_k^{(\al,\be)}}
\newcommand{\paaz}{P_0^{(\al,\al)}}
\newcommand{\paakp}{P_{k+1}^{(\al,\al)}}
\newcommand{\paakm}{P_{k-1}^{(\al,\al)}}

\newcommand{\pookm}{P_{k-1}^{(\al+1,\al+1)}}
\newcommand{\pook}{P_{k}^{(\al+1,\al+1)}}
\newcommand{\pozj}{P_{j}^{(\al+1,\al)}}
\newcommand{\pozk}{P_{k}^{(\al+1,\al)}}
\newcommand{\pozkm}{P_{k-1}^{(\al+1,\al)}}
\newcommand{\pzok}{P_{k}^{(\al,\al+1)}}

\newcommand{\tab}{t^{(\al,\be)}}

\newcommand{\too}{t^{(\al+1,\al+1)}}
\newcommand{\tookmkm}{\too_{k-1,k-1}}

\newcommand{\toojmjm}{\too_{j-1,j-1}}
\newcommand{\toojmim}{\too_{j-1,i-1}}
\newcommand{\toojmi}{\too_{j-1,i}}
\newcommand{\tozjj}{t^{(\al+1,\al)}_{j,j}}
\newcommand{\tozkmkm}{t^{(\al+1,\al)}_{k-1,k-1}}

\newcommand{\tookmi}{t^{(\al+1,\al+1)}_{k-1,i}}
\newcommand{\tookmmim}{t^{(\al+1,\al+1)}_{k-2,i-1}}
\newcommand{\tookmo}{t^{(\al+1,\al+1)}_{k-1,1}}

\begin{document}

\begin{abstract}
The Lov\'asz theta function provides a lower bound for the chromatic
number of finite graphs based on the solution of a semidefinite
program. In this paper we generalize it so that it gives a lower bound
for the measurable chromatic number of distance graphs on compact metric spaces.

In particular we consider distance graphs on the unit sphere. There we
transform the original infinite semidefinite program into an infinite
linear program which then turns out to be an extremal question about
Jacobi polynomials which we solve explicitly in the limit. As an
application we derive new lower bounds for the measurable chromatic
number of the Euclidean space in dimensions $10, \dots, 24$ and we
give a new proof that it grows exponentially with the dimension.
\end{abstract}

\maketitle

\markboth{C.~Bachoc, G.~Nebe, F.M.~de Oliveira Filho,
F.~Vallentin}{Lower bounds for measurable chromatic numbers}

\section{Introduction}

The \defi{chromatic number of the $n$-dimensional Euclidean space} is
the minimum number of colors needed to color each point of $\R^n$ in
such a way that points at distance~$1$ from each other receive
different colors. It is the chromatic number of the graph with vertex
set $\R^n$ and in which two vertices are adjacent if their distance
is~$1$. We denote it by $\chi(\R^n)$.

A famous open question is to determine the chromatic number of the
plane. In this case, it is only known that $4 \leq \chi(\R^2) \leq 7$,
where lower and upper bounds come from simple geometric
constructions. In this form the problem was considered, e.g., by
Nelson, Isbell, Erd\H{o}s, and Hadwiger. For historical
remarks and for the best known bounds in other dimensions we refer to
Sz\'ekely's survey article \cite{Szek}. The first exponential
asymptotic lower bound is due to Frankl and
Wilson~\cite[Theorem 3]{FW}. Currently the best known asymptotic
lower bound is due to Raigorodskii~\cite{R} and the best known
asymptotic upper bound is due to Larman and
Rogers~\cite{LR}:
\begin{equation*}
(1.239\ldots + o(1))^n \leq \chi(\R^n) \leq (3 + o(1))^n.
\end{equation*}

In this paper we study a variant of the chromatic number of $\R^n$,
namely the measurable chromatic number. The \defi{measurable chromatic
number} of~$\R^n$ is the smallest number $m$ such that $\R^n$ can be
partitioned into $m$ Lebesgue measurable stable sets. Here we call a
set $C \subseteq \R^n$ \defi{stable} if no two points in $C$ lie at
distance~$1$ from each other. In other words, we impose that the sets
of points having the same color have to be measurable. We denote the
measurable chromatic number of $\R^n$ by $\chim(\R^n)$. One reason to
study the measurable chromatic number is that then stronger analytic
tools are available.

The study of the measurable chromatic number started with
Falconer~\cite{Fal}, who proved that $\chim(\R^2) \geq 5$.
The measurable chromatic number is at least the chromatic number and
it is amusing to notice that in case of strict inequality the
construction of an optimal coloring necessarily uses the axiom of
choice.

\smallskip

Related to the chromatic number of the Euclidean space is the
chromatic number of the unit sphere $S^{n-1} = \{x \in \R^n : x \cdot
x = 1\}$. For $-1 < t < 1$, we consider the graph $G(n, t)$ whose
vertices are the points of $S^{n-1}$ and in which two points are
adjacent if their inner product $x \cdot y$ equals $t$. The chromatic
number of $G(n,t)$ and its measurable version, denoted by
$\chi(G(n,t))$ and $\chim(G(n,t))$ respectively, are defined as in
the Euclidean case.

The chromatic number of this graph was studied by
Lov\'asz~\cite{Lov2}, in particular in the case when $t$ is
small. He showed that
\begin{equation*}
\begin{split}
&  n \leq \chi(G(n,t))\quad \text{ for $-1 < t < 1$,}\\
& \chi(G(n,t)) \leq n+1\quad \text{ for $-1 < t \leq -1/n$}.
\end{split}
\end{equation*}
Frankl and Wilson~\cite[Theorem 6]{FW} showed that
\begin{equation*}
(1+o(1))(1.13)^n \leq \chim(G(n,0)) \leq 2^{n-1}.
\end{equation*}

The (measurable) chromatic number of $G(n,t)$ provides a lower bound
for the one of $\R^n$: After appropriate scaling, every proper
coloring of $\R^n$ intersected with the unit sphere $S^{n-1}$ gives a
proper coloring of the graph $G(n,t)$, and measurability is preserved
by the intersection.

\smallskip

In this paper we present a lower bound for the measurable chromatic
number of $G(n,t)$. As an application we derive new lower bounds for
the measurable chromatic number of the Euclidean space in dimensions
$10,\ldots,24$ and we give a new proof that it grows exponentially
with the dimension.

The lower bound is based on a generalization of the Lov\'asz theta
function (Lov\'asz~\cite{Lov1}), which gives an upper bound to the
stability number of a finite graph. Here we aim at generalizing the
theta function to \defi{distance graphs} in compact metric
spaces. These are graphs defined on all points of the metric space
where the adjacency relation only depends on the distance.

The paper is structured as follows: In Section~\ref{sec:fractional} we
define the stability number and the fractional measurable chromatic
number and give a basic inequality involving them. Then, after
reviewing Lov\'asz' original formulation of the theta function in
Section~\ref{sec:original theta}, we give our generalization in
Section~\ref{sec:generalization}.  Like the original theta function
for finite graphs, it gives an upper bound for the stability
number. Moreover, in the case of the unit sphere, it can be explicitly
computed, thanks to classical results on spherical harmonics. The
material needed for spherical harmonics is given in
Section~\ref{sec:spherical harmonics} and an explicit formulation for
the theta function of $G(n,t)$ is given in
Section~\ref{sec:generalization2}.

In Section~\ref{sec:analytic} we choose specific values of $t$ for
which we can analytically compute the theta function of $G(n,t)$.
This allows us to compute the limit of the theta function for the
graph $G(n,t)$ as $t$ goes to $1$ in Section~\ref{sec:new lower
bounds}. This gives improvements on the best known lower bounds for
$\chim(\R^n)$ in several dimensions. Furthermore this gives a new
proof of the fact that $\chim(\R^n)$ grows exponentially with
$n$. Although this is an immediate consequence of the result of
Frankl and Wilson (and of Raigorodskii, and also of a
result of Frankl and R\"odl \cite{FR}) and our bound of
$1.165^n$ is not an improvement, our result is an easy consequence of
the methods we present. Moreover, we think that our proof is of
interest because the methods we use here are radically different from
those used before. In particular, they can be applied to other metric
spaces.

In Section~\ref{sec:further} we point out how to apply our
generalization to distance graphs in other compact metric spaces,
endowed with the continuous action of a compact group. Finally in
Section~\ref{sec:second} we conclude by showing the relation between
our generalization of the theta function and the theta function for
finite graphs of $G(n,t)$ and by showing the relation between our
generalization and the linear programming bound for spherical codes
established by Delsarte, Goethals, and Seidel~\cite{DGS}.

\section{The fractional chromatic number and the stability number}
\label{sec:fractional}

Let $G = (V, E)$ be a finite or infinite graph whose vertex set is
equipped with the measure $\mu$. We assume that the measure of $V$ is
finite. In this section we define the stability number and the
measurable fractional chromatic number of $G$ and derive the basic
inequality between these two invariants. In the case of a finite graph
one recovers the classical notions if one uses the uniform measure
$\mu(C) = |C|$ for $C \subseteq V$.

Let $L^2(V)$ be the Hilbert space of real-valued square-integrable
functions defined over $V$ with inner product
\begin{equation*}
(f, g) = \int_V f(x) g(x)\, d\mu(x)
\end{equation*}
for $f, g \in L^2(V)$. The constant function $1$ is measurable and its
squared norm is the number $(1,1) = \mu(V)$. The characteristic
function of a subset $C$ of $V$ we denote by $\chi^C\colon V \to \{ 0,
1\}$.

A subset $C$ of $V$ is called a \defi{measurable stable set} if $C$ is
a measurable set and if no two vertices in $C$ are adjacent. The
\defi{stability number} of $G$ is
\begin{equation*}
\alpha(G) = \sup \{ \mu(C) : C \subseteq V \text{ is a measurable
stable set} \}.
\end{equation*}
Similar measure-theoretical notions of the stability number have been
considered before by other authors for the case in which $V$ is the
Euclidean space $\R^n$ or the sphere $S^{n-1}$. We refer the reader to
the survey paper of Sz\'{e}kely~\cite{Szek} for more information and
further references.

The \defi{fractional measurable chromatic number} of $G$ is denoted by
$\chim^*(G)$. It is the infimum of $\lambda_1 + \dots + \lambda_k$
where $k \geq 0$ and $\lambda_1, \dots, \lambda_k$ are nonnegative
real numbers such that there exist measurable stable sets $C_1, \dots,
C_k$ satisfying
\begin{equation*}
\lambda_1 \chi^{C_1} + \cdots + \lambda_k \chi^{C_k} = 1.
\end{equation*}
Note that the measurable fractional chromatic number of the graph $G$
is a lower bound for its measurable chromatic number.

\begin{proposition}
We have the following basic inequality between the stability number
and the measurable fractional chromatic number of a graph $G = (V,E)$:
\begin{equation}
\label{eq:chi-alpha}
\alpha(G) \chim^*(G) \geq \mu(V).
\end{equation}
So, any upper bound for $\alpha(G)$ provides a lower bound for
$\chim^*(G)$.
\end{proposition}

\begin{proof}
Let $\lambda_1, \dots, \lambda_k$ be nonnegative real numbers and
$C_1, \dots, C_k$ be measurable stable sets such that $\lambda_1
\chi^{C_1} + \cdots + \lambda_k \chi^{C_k} = 1$. Since $C_i$ is
measurable, its characteristic function $\chi^{C_i}$ lies in $L^2(V)$.
Hence
\begin{equation*}
\begin{split}
(\lambda_1 + \cdots + \lambda_k) \alpha(G) & \geq
\lambda_1 \mu(C_1) + \cdots + \lambda_k \mu(C_k)\\
 &= \lambda_1 (\chi^{C_1}, 1) + \cdots + \lambda_k (\chi^{C_k},
1)\\
&= (1, 1)\\
& = \mu(V).
\qedhere
\end{split}
\end{equation*}
\end{proof}

\section{The Lov\'asz theta function for finite graphs}
\label{sec:original theta}

In the celebrated paper \cite{Lov1}  Lov\'asz introduced the theta
function for finite graphs. It is an upper bound for the stability
number which one can efficiently compute using semidefinite
programming. In this section we review its definition and properties,
which we generalize in Section~\ref{sec:generalization}. 

The \defi{theta function} of a  graph $G = (V, E)$
is defined by
\begin{equation}
\label{original theta}
\begin{split}
  \vartheta(G) = \max\big\{ & \sum\nolimits_{x \in V}\sum\nolimits_{y \in V} K(x,y) :  \\
& \quad \text{$K \in \R^{V \times V}$ is positive semidefinite},\\
& \quad \sum_{x\in V} K(x,x) = 1,\\
& \quad \text{$K(x,y) = 0$ if $\{x,y\} \in E$}\big\}.
\end{split}
\end{equation}

\begin{theorem}
\label{th:original}
For any finite graph $G$, $\vartheta(G) \geq \alpha(G)$.
\end{theorem}

Although this result follows from \cite[Lemma 3]{Lov1} and
\cite[Theorem 4]{Lov1}, we give a proof here to stress the analogy
between the finite case and the more general case we consider in our
generalization Theorem~\ref{thm:theta}.

\begin{proof}[Proof of Theorem~\textup{\ref{th:original}}]
Let $C \subseteq V$ be a stable set. Consider the characteristic
function $\chi^C\colon V \to \{0,1\}$ of $C$ and define the matrix $K
\in \R^{V \times V}$ by
\begin{equation*}
K(x,y) = \frac{1}{|C|} \chi^C(x) \chi^C(y).
\end{equation*}
Notice $K$ satisfies the conditions in \eqref{original theta}. 
Moreover, we have $\sum_{x \in
V}\sum_{y \in V} K(x,y) = |C|$, and so $\vartheta(G) \geq |C|$.
\end{proof}

\begin{remark}
\label{rem:alt def}
There are many equivalent definitions of the theta
function. Possible alternatives are reviewed by Knuth in
\cite{K}. We use the one of \cite[Theorem 4]{Lov1}.
\end{remark}

If the graph $G$ has a nontrivial automorphism group, it is not
difficult to see that one can restrict oneself in \eqref{original
theta} to the functions $K$ which are invariant under the action of
any subgroup $\Gamma$ of $\Aut(G)$, where
$\Aut(G)$ is the \defi{automorphism group} of $G$, i.e., it is the
group of all permutations of $V$ that preserve adjacency.  Here we say
that $K$ is \defi{invariant under $\Gamma$} if $K(\gamma x, \gamma y)
= K(x,y)$ holds for all $\gamma \in \Gamma$ and all $x,y \in V$.  If
moreover $\Gamma$ acts transitively on $G$, the second condition
$\sum_{x\in V} K(x,x) = 1$ is equivalent to $K(x,x)=1/|V|$ for all $x\in
V$.

\section{A generalization of the Lov\'asz theta function for distance
  graphs on compact metric spaces}
\label{sec:generalization}

We assume that $V$ is a compact metric space with distance function
$d$. We moreover assume that $V$ is equipped with a nonnegative, Borel regular measure
$\mu$ for which $\mu(V)$ is finite. Let $D$ be a closed subset of the
image of $d$.  We define the graph $G(V,D)$ to be the graph with
vertex set $V$ and edge set $E = \{\{x,y\} : d(x,y)\in D\}$. 

The elements of the space $\Ccal(V \times V)$ consisting of all
continuous functions $K : V \times V \to \R$ are called
\defi{continuous Hilbert-Schmidt kernels}; or \defi{kernels} for short. In
the following we only consider \defi{symmetric} kernels, i.e., kernels
$K$ with $K(x,y) = K(y,x)$ for all $x, y \in V$.  A kernel $K \in \Ccal(V
\times V)$ is called \defi{positive} if, for any nonnegative integer
$m$, any points $x_1, \ldots, x_m \in V$, and any real numbers
$u_1, \ldots, u_m$, we have
\begin{equation*}
\sum_{i = 1}^m \sum_{j = 1}^m K(x_i,x_j) u_i u_j \geq 0.
\end{equation*}

We are now ready to extend the definition \eqref{original theta} of the Lov\'asz theta
function to the graph $G(V,D)$. We define
\begin{equation}
\label{eq:new theta}
\begin{split}
  \vartheta(G(V,D)) = \sup\big\{ & \int_V\int_V  K(x,y)\, d\mu(x)d\mu(y):  \\
& \quad \text{$K \in \Ccal(V \times V)$ is positive},\\
& \quad \int_V K(x,x)\, d\mu(x)=1,\\ 
& \quad \text{$K(x,y) = 0$ if $d(x,y)\in D$}\big\}.
\end{split}
\end{equation}

\begin{theorem}
\label{thm:theta}
The theta function is an upper bound for the stability number, i.e.,
\begin{equation*}
\vartheta(G(V,D)) \geq \alpha(G(V,D)).
\end{equation*}
\end{theorem}

\begin{proof} 
Fix $\varepsilon > 0$ arbitrarily. Let $C \subseteq V$ be a stable set
such that $\mu(C) \geq \alpha(G(V, D)) - \varepsilon$. Since $\mu$ is
regular, we may assume that $C$ is closed, as otherwise we could
find a stable set with measure closer to $\alpha(G(V, D))$ and use a
suitable inner-approximation of it by a closed set.

Note that, since $C$ is compact and stable, there must exist a number
$\beta > 0$ such that $| d(x, y) - \delta | > \beta$ for all $x, y \in
C$ and $\delta \in D$. But then, for small enough $\xi > 0$, the set
\begin{equation*}
B(C, \xi) = \{ x \in V : d(x, C) < \xi \},
\end{equation*}
where $d(x, C)$ is the distance from $x$ to the closed set $C$, is
stable. Moreover, notice that $B(C, \xi)$ is open and that, since it
is stable, $\mu(B(C, \xi)) \leq \alpha(G(V, D))$.

Now, the function $f\colon V \to [0, 1]$ given by
\begin{equation*}
f(x) = \xi^{-1} \cdot \max \{ \xi - d(x, C), 0 \}
\end{equation*}
for all $x \in V$ is continuous and such that $f(C) = 1$ and
$f(V \setminus B(C, \xi)) = 0$.
So the kernel $K$ given by
\begin{equation*}
K(x, y) = \frac{1}{(f, f)} f(x) f(y)
\end{equation*}
for all $x, y \in V$ is feasible in~\eqref{eq:new theta}. 

Let us estimate the objective value of $K$. Since we have
\begin{equation*}
(f, f) \leq \mu(B(C, \xi)) \leq \alpha(G(V,  D))
\end{equation*}
and
\begin{equation*}
\int_V \int_V f(x) f(y)\, d\mu(x) d\mu(y) \geq \mu(C)^2 \geq
(\alpha(G(V, D)) - \varepsilon)^2,
\end{equation*}
we finally have
\begin{equation*}
\int_V \int_V K(x, y)\, d\mu(x) d\mu(y) \geq \frac{(\alpha(G(V, D)) -
  \varepsilon)^2}{\alpha(G(V, D))}
\end{equation*}
and, since $\varepsilon$ is arbitrary, the theorem follows.
\end{proof}

Let us now assume that a compact group $\Gamma$ acts continuously on
$V$, preserving the distance $d$. Then, if $K$ is a feasible solution
for \eqref{eq:new theta}, so is $(x,y)\mapsto K(\gamma x,\gamma y)$
for all $\gamma \in \Gamma$. Averaging on $\Gamma$ leads to a
$\Gamma$-invariant feasible solution
\begin{equation*}
\overline{K}(x,y)=\int_{\Gamma} K(\gamma x,\gamma y)\, d\gamma,
\end{equation*}
where $d\gamma$ denotes the Haar measure on $\Gamma$ normalized so
that $\Gamma$ has volume $1$. Moreover, observe that the objective
value of $\overline{K}$ is the same as that of $K$. Hence we can
restrict ourselves in \eqref{eq:new theta} to $\Gamma$-invariant
kernels.  If moreover $V$ is  homogeneous under the action of $\Gamma$, the second
condition in~\eqref{eq:new theta} may be replaced by $K(x,x)=1/\mu(V)$ for all $x\in V$. 

We are mostly interested in the case in which $V$ is the unit sphere
$S^{n-1}$ endowed with the Euclidean metric of $\R^n$, and in which $D$ is a
singleton. If $D=\{\delta\}$ and $\delta^2=2-2t$, so that
$d(x,y)=\delta$ if and only if $x\cdot y=t$, the graph $G(S^{n-1}, D)$
is denoted by $G(n,t)$. Since the unit sphere is homogeneous under the
action of the orthogonal group $\ort(\R^n)$, the previous remarks
apply.

\section{Harmonic analysis on the unit sphere}
\label{sec:spherical harmonics}

It turns out that the continuous positive Hilbert-Schmidt kernels on
the sphere have a nice description coming from classical results of
harmonic analysis reviewed in this section. This allows for the
calculation of $\vartheta(G(n,t))$. For information on spherical
harmonics we refer to \cite[Chapter 9]{AAR} and \cite{VK}.

The unit sphere $S^{n-1}$ is homogeneous under the action of the
orthogonal group $\ort(\R^n) = \{A \in \R^{n \times n} : A^t A =
I_n\}$, where $I_n$ denotes the identity matrix. Moreover, it is
two-point homogeneous, meaning that the orbits of $\ort(\R^n)$ on
pairs of points are characterized by the value of their inner
product. The orthogonal group acts on $L^2(S^{n-1})$ by $Af(x) =
f(A^{-1}x)$, and $L^2(S^{n-1})$ is equipped with the standard
$\ort(\R^n)$-invariant inner product
\begin{equation}
\label{eq:invariant inner product}
(f,g) = \int_{S^{n-1}} f(x)g(x)\, d\omega(x)
\end{equation}
for the standard surface measure $\omega$. The surface area of the
unit sphere is $\omega_n = (1, 1) = 2\pi^{n/2} / \Gamma(n/2)$.

It is a well-known fact (see e.g.\ \cite[Chapter
9.2]{VK}) that the Hilbert space $L^2(S^{n-1})$ decomposes under the
action of $\ort(\R^n)$ into orthogonal subspaces
\begin{equation}
\label{eq:dec}
L^2(S^{n-1}) = H_0\perp H_1\perp H_2 \perp \ldots,
\end{equation}
where $H_k$ is isomorphic to the $\ort(\R^n)$-irreducible space
\begin{equation*}
\Harm_k = \Big\{f \in \R[x_1, \ldots, x_n] : \text{$f$
homogeneous}, \deg f = k, \sum_{i=1}^n \frac{\partial^2}{\partial
x_i^2} f = 0\Big\}
\end{equation*}
of harmonic polynomials in $n$ variables which are homogeneous and
have degree $k$.  We set $h_k=\dim(\Harm_k) = \binom{n+k-1}{n-1} -
\binom{n+k-3}{n-1}$. The equality in \eqref{eq:dec} means that every
$f \in L^2(S^{n-1})$ can be uniquely written in the form $f =
\sum_{k = 0}^\infty p_k$, where $p_k \in H_k$, and where the
convergence is in the $L^2$-norm.

The \defi{addition formula} (see e.g.~\cite[Chapter 9.6]{AAR}) plays a
central role in the characterization of $\ort(\R^n)$-invariant
kernels: For any orthonormal basis $e_{k,1},\ldots,e_{k,{h_k}}$ of
$H_k$ and for any pair of points $x,y \in S^{n-1}$ we have
\begin{equation}
\label{add form}
\sum_{i=1}^{h_k}
e_{k,i}(x)e_{k,i}(y) = \frac{h_k}{\omega_n} 
\paak(x \cdot y),
\end{equation}
where $\paak$ is the normalized Jacobi polynomial of degree $k$ with
parameters $(\al,\al)$, with $\paak(1) = 1$ and $\al = (n-3)/2$. The
\defi{Jacobi polynomials} with parameters $(\alpha, \beta)$ are
orthogonal polynomials for the weight function $(1-u)^\al(1+u)^\be$ on
the interval $[-1,1]$. We denote by $\pabk$ the normalized Jacobi
polynomial of degree $k$  with normalization $\pabk(1) = 1$. 

In \cite[Theorem 1]{Sch} Schoenberg gave a characterization of
the continuous kernels which are positive and $\ort(\R^n)$-invariant:
They are those which lie in the cone spanned by the kernels $(x,y)
\mapsto \paak(x \cdot y)$. More precisely, a continuous kernel $K \in
\Ccal(S^{n-1} \times S^{n-1})$ is $\ort(\R^n)$-invariant and positive
if and only if there exist nonnegative real numbers $f_0, f_1, \dots$
such that $K$ can be written as
\begin{equation}
\label{eq:schoenberg}
K(x, y) = \sum_{k = 0}^\infty f_k \paak(x \cdot y),
\end{equation}
where the convergence is absolute and uniform.

\section{The theta function of $G(n,t)$}
\label{sec:generalization2}

We obtain from Section \ref{sec:generalization} in the case
$V=S^{n-1}$, $D=\{\sqrt{2-2t}\}$, and $\Gamma = \ort(\R^n)$, the
following characterization of the theta function of the graph $G(n,t)
= G(S^{n-1},D)$:
\begin{equation}
\label{eq:theta G(n,t)}
\begin{split}
  \vartheta(G(n,t)) = \max\big\{ & \intsn K(x,y)\, d\omega(x)d\omega(y):  \\
& \quad \text{$K \in \Ccal(S^{n-1} \times S^{n-1})$ is positive},\\
& \quad \text{$K$ is invariant under $\ort(\R^n)$},\\
& \quad \text{$K(x,x) = 1/\omega_n$ for all $x \in S^{n-1}$,}\\
& \quad \text{$K(x,y) = 0$ if $x \cdot y = t$}\big\}.
\end{split}
\end{equation}
(It will be clear later that the maximum above indeed exists.)

\begin{corollary}
\label{cor:lowerboundchi}
We have
\begin{equation*}
\omega_n / \vartheta(G(n, t)) \leq \chim^*(G(n, t)).
\end{equation*}
\end{corollary}

\begin{proof}
Immediate from Theorem \ref{thm:theta} and the considerations in
Section~\ref{sec:fractional}.
\end{proof}

A result of de Bruijn and Erd\H{o}s~\cite{BE} implies that
the chromatic number of $G(n,t)$ is attained by a finite 
subgraph of it. So one might wonder if computing the theta
function for a finite subgraph of $G(n,t)$ could give a better
bound than the previous corollary. This is not the case as we will
show in Section~\ref{sec:second}.

The theta function for finite graphs has the important property that
it can be computed in polynomial time, in the sense that it can be
approximated with arbitrary precision using semidefinite
programming. We now turn to the problem of computing the
generalization~\eqref{eq:theta G(n,t)}.

First, we apply Schoenberg's characterization~\eqref{eq:schoenberg} of
the continuous kernels which are $\ort(\R^n)$-invariant and positive.
This transforms the original formulation~\eqref{eq:new theta}, which
is a semidefinite programming problem in infinitely many variables
having infinitely many constraints, into the following linear
programming problem with optimization variables $f_k$:
\begin{equation}
\label{eq:primal}
\begin{split}
  \vartheta(G(n,t)) = \max\big\{ & \omega_n^2 f_0:  \\
& \quad \text{$f_k \geq 0$ for $k = 0, 1, \ldots$},\\
& \quad \text{$\sum\limits_{k = 0}^\infty f_k = 1/\omega_n$,}\\
& \quad \text{$f_0 + \sum\limits_{k = 1}^\infty f_k \paak(t) = 0$}\big\},
\end{split}
\end{equation}
where $\alpha = (n-3)/2$.  

To obtain \eqref{eq:primal} we simplified the objective function in
the following way.  Because of the orthogonal
decomposition~\eqref{eq:dec} and because the subspace $H_0$ contains
only the constant functions, we have
\begin{equation*}
\intsn \sum\limits_{k = 0}^\infty f_k \paak(x \cdot y)\, d\omega(x) d\omega(y) = \omega_n^2 f_0.
\end{equation*}
We furthermore used $\paaz = 1$ and $\paak(1) = 1$.

\begin{theorem}
\label{thm:optimal solution}
Let $m(t)$ be the minimum of $\paak(t)$ for $k = 0, 1, \ldots$ Then the
optimal value of~\eqref{eq:primal} is equal to
\begin{equation*}
\vartheta(G(n,t)) = \omega_n \frac{m(t)}{m(t) - 1}.
\end{equation*}
\end{theorem}

\begin{proof}
  We first claim that the minimum $m(t)$ exists and is
  negative. Indeed, if $\paak(t) \geq 0$ for all $k \geq 1$,
  then~\eqref{eq:primal} either has no solution (in the case that all
  $\paak(t)$ are positive) or $f_0 = 0$ in any solution, which
  contradicts Theorem~\ref{thm:theta}.  So we know that for some $k
  \geq 1$, $\paak(t) < 0$. This, combined with the fact that $\paak(t)$
  goes to zero as $k$ goes to infinity (cf.~\cite[Chapter 6.6]{AAR}
  or~\cite[Chapter 8.22]{Szeg}), proves the claim.

  Let $k^*$ be so that $m(t) = P^{(\alpha, \alpha)}_{k^*}(t)$. It is
  easy to see that there is an optimal solution of~\eqref{eq:primal} in which only
  $f_0$ and $f_{k^*}$ are positive. Hence, solving the resulting system
\begin{eqnarray*}
f_0 + f_{k^*} & = & 1/\omega_n\\
f_0 + f_{k^*} m(t) & = &  0
\end{eqnarray*}
gives $f_0 = m(t)/(\omega_n(m(t)-1))$ and
$f_{k^*} = -1/(\omega_n (m(t)-1))$ and the theorem follows.
\end{proof}

\begin{example}
The minimum of $\paak(0.9999)$ for $\alpha = (24-3)/2$ is attained at
$k = 1131$. It is a rational number and its first decimal digits are $-0.00059623$.
\end{example}

\section{Analytic solutions}
\label{sec:analytic}

In this section we compute the value
\begin{equation*}
m(t) = \min\{\paak(t) : k = 0, 1, \ldots\}
\end{equation*}
for specific values of $t$. Namely we choose $t$ to be the largest
zero of an appropriate Jacobi polynomial.

Key for the discussion to follow is the \defi{interlacing property} of
 the zeroes of orthogonal polynomials. It says (cf.~\cite[Theorem
 3.3.2]{Szeg}) that between any pair of consecutive zeroes of $\paak$
 there is exactly one zero of $\paakm$.

We denote the zeros of $\pabk$ by $\tab_{k,j}$ with $j = 1, \ldots, k$
and with the increasing ordering $\tab_{k,j} < \tab_{k,j+1}$. We shall
need the following collection of identities:
\begin{equation}
\label{e1}
(1-u^2)\frac{d^2\paak}{du^2} - (2\al+2)u\frac{d\paak}{du} + k(k+2\al+1)\paak=0,
\end{equation}
\begin{align}
\label{e2a}
(-1)^k \paak(-u)& =\paak(u),\\
\label{e2b}
(-1)^k (\al+1)\pzok(-u) & = (k + \al + 1) \pozk(u),\\
\label{e3}
(2\al+2)\frac{d\paak}{du} &= k(k+2\al+1) \pookm, \\
\label{e4}
(2\al+2)\pzok &= (k+2\al+2) \pook-k \pookm,\\
\label{e5}
(2k+2\al+2)\pozk &= (k+2\al+2) \pook+k \pookm,\\
\label{e6}
(k+\al+1)\pozk&= (\al+1)\frac{\paak -\paakp}{1-u}.
\end{align}
They can all be found in \cite[Chapter 6]{AAR}, although with
different normalization. Formula \eqref{e1} is \cite[(6.3.9)]{AAR};
\eqref{e2a} and \eqref{e2b} are \cite[(6.4.23)]{AAR}; \eqref{e3} is
\cite[(6.3.8)]{AAR}, \eqref{e4} is \cite[(6.4.21)]{AAR}; \eqref{e5}
follows by the change of variables $u \mapsto -u$ from \eqref{e4} and
\eqref{e2a}, \eqref{e2b}; \eqref{e6} is \cite[(6.4.20)]{AAR}.

\begin{proposition}
\label{prop:tookmkm}
Let $t = \tookmkm$ be the largest zero of the Jacobi polynomial
$\pookm$. Then, $m(t) = \paak(t)$.
\end{proposition}

\begin{proof} 
We start with the following crucial observation: From \eqref{e3}, $t$
is a zero of the derivative of $\paak$. Hence it is a minimum of
$\paak$ because it is the last extremal value in the interval $[-1,1]$ and
because $\pook(1) = 1$, whence (using~\eqref{e3}) $\paak(u)$ is
increasing on $[t, 1]$.

Now we prove that $\paak(t) < \paaj(t)$ for all $j\neq k$ where we
treat the cases $j<k$ and $j>k$ separately.

It turns out that the sequence $\paaj(t)$ is decreasing for $j\leq
k$. From \eqref{e6}, the sign of $\paaj(t)-\paajp(t)$ equals the
sign of $\pozj(t)$. We have the inequalities
\begin{equation*}
\tozjj \leq \tozkmkm < \tookmkm = t.
\end{equation*}
The first one is a consequence of the interlacing property. From
\eqref{e5} one can deduce that $\pozkm$ has exactly one zero in the
interval $[\tookmmim, \tookmi]$ since it changes sign at the extreme
points of it, and by the same argument $\pozkm$ has a zero left to
$\tookmo$. Thus, $\tozkmkm < \tookmkm = t$. So $t$ lies to the right of the
largest zero of $\pozj$ and hence $\pozj(t) > 0$ which shows that
$\paaj(t) - \paajp(t)>0$ for $j < k$.

Let us consider the case $j>k$. The inequality \cite[(6.4.19)]{AAR}
implies that
\begin{equation}\label{ee2}
\text{for all $j > k$, $\quad\paak(\tookmkm) < \paaj(\toojmjm)$}.
\end{equation}
The next observation, which finishes the proof of the lemma, is
stated in \cite[(6.4.24)]{AAR} only for the case $\alpha = 0$:
\begin{equation}\label{ee1}
\text{for all $j\geq 2$, $\quad\min\{\paaj(u) : u \in [0,1]\} = \paaj(\toojmjm)$}.
\end{equation}
To prove it consider
\begin{equation*}
g(u)=\paaj(u)^2 + \frac{1-u^2}{j(j+2\al+1)}\Big(\frac{d\paaj}{du}\Big)^2.
\end{equation*}
Applying \eqref{e1} in the computation of $g'$ shows that 
\begin{equation*}
g'(u)= \frac{(4\al+2)u}{j(j+2\al+1)}\Big(\frac{d\paaj}{du}\Big)^2.
\end{equation*}
The polynomial $g'$ takes positive values on $[0,1]$ and hence $g$ is
increasing on this interval. In particular,
\begin{equation*}
\text{$g(\toojmim) < g(\toojmi)\quad$ for all $i\leq j-1$ with $\toojmim \geq 0$},
\end{equation*}
which simplifies to
\begin{equation*}
  \text{$\paaj(\toojmim)^2 < \paaj(\toojmi)^2\quad$}.
\end{equation*}
Since $\toojmi$ are the local extrema of $\paaj$, we have
proved \eqref{ee1}.
\end{proof}

\section{New lower bounds for the Euclidean space}
\label{sec:new lower bounds}

In this section we give new lower bounds for the measurable chromatic
number of the Euclidean space for dimensions $10, \ldots, 24$. This
improves on the previous best known lower bounds due to Sz\'ekely
and Wormald \cite{SW}. Table~\ref{table:thetable} compares the
values.  Furthermore we give a new proof that the measurable chromatic
number grows exponentially with the dimension.

For this we give a closed expression for $\lim_{t \to 1} m(t)$ which
involves the Bessel function $J_\al$ of the first kind of order
$\alpha = (n-3)/2$ (see e.g.\ \cite[Chapter 4]{AAR}). The appearance
of Bessel functions here is due to the fact that the largest zero of
the Jacobi polynomial $\paak$ behaves like the first positive zero
$j_\al$ of the Bessel function $J_\al$. More precisely, it is
known~\cite[Theorem 4.14.1]{AAR} that, for the largest zero $t^{(\alpha+1,\beta)}_{k,k} = 
\cos\theta_k$ of the polynomial $P^{(\alpha+1,\beta)}_{k}$,
\begin{equation}\label{j1}
\lim_{k\to\infty} k \theta_k = j_{\al+1}
\end{equation}
and, with our normalization (cf.\ \cite[Theorem 4.11.6]{AAR}),
\begin{equation}\label{j2}
\lim_{k\to\infty} \paak\left(\cos\frac{u}{k}\right) = 2^\al\Gamma(\al+1)\frac{J_\al(u)}{u^\al}.
\end{equation}

\begin{theorem}
We have
\begin{equation*}
\lim_{t \to 1} m(t) = 2^\al \Gamma(\al+1)
\frac{J_{\al}(j_{\al+1})}{(j_{\al+1})^\al}.
\end{equation*}

\end{theorem}

\begin{proof}
First we show that
\begin{equation}
\label{eq:limit}
\lim_{k \to \infty} P^{(\alpha,\alpha)}_k(t^{(\alpha+1,\beta)}_{k-1,k-1}) = 2^{\alpha}\Gamma(\alpha+1)\frac{J_{\alpha}(j_{\alpha+1})}{(j_{\alpha+1})^{\alpha}}.
\end{equation}
We estimate the difference 
\begin{equation*}
|\paak(t^{(\alpha+1,\beta)}_{k-1,k-1}) - 2^\al
 \Gamma(\al+1)\frac{J_{\al}(j_{\al+1})}{(j_{\al+1})^\al}|,
\end{equation*}
that we upper bound by 
\begin{equation*}
\begin{split}
&
\left|\paak(t^{(\alpha+1,\beta)}_{k-1,k-1})-\paak\left(\cos\frac{j_{\al+1}}{k}\right)\right|\\
& \quad +\left|\paak\left(\cos\frac{j_{\al+1}}{k}\right)- 2^\al
 \Gamma(\al+1)\frac{J_{\al}(j_{\al+1})}{(j_{\al+1})^\al}\right|.
\end{split}
\end{equation*}

The second term tends to $0$ from \eqref{j2}. Define $\theta_{k-1}$ by $t^{(\alpha+1,\beta)}_{k-1,k-1}=\cos
\theta_{k-1}$. By the mean value theorem we have
\begin{align*}
& \left|\paak(t^{(\alpha+1,\beta)}_{k-1,k-1})-\paak\left(\cos\frac{j_{\al+1}}{k}\right)\right|\\
\leq & \quad \Big(\max_{u \in [-1,1]}\big|\frac{d\paak}{du}\big|\Big) \big|\cos \theta_{k-1} - \cos\frac{j_{\al+1}}{k}\big|\\
\leq & \quad \Big(\max_{u \in [-1,1]}\big|\frac{d\paak}{du}\big|\Big) \big(\max_{\theta \in I_k}|\sin \theta|\big) \big|\theta_{k-1} - \frac{j_{\al+1}}{k}\big|,
\end{align*}
where $I_k$ denotes the interval with extremes $\theta_{k-1}$ and
$\frac{j_{\al+1}}{k}$. Then, with \eqref{j1},
\begin{align*}
\theta_{k-1} - \frac{j_{\al+1}}{k}&=\theta_{k-1} -
\frac{j_{\al+1}}{k-1} +\frac{j_{\al+1}}{k(k-1)}\\
&=\frac{1}{k-1}((k-1)\theta_{k-1} -
j_{\al+1})+\frac{j_{\al+1}}{k(k-1)}=o\left(\frac{1}{k}\right),
\end{align*}
and for all $\theta\in I_k$
\begin{equation*}
|\sin \theta| \leq |\theta|\leq
\frac{j_{\al+1}}{k} + \big|\theta_{k-1} - \frac{j_{\al+1}}{k}\big| = O\left(\frac{1}{k}\right).
\end{equation*}
From \eqref{e3}, 
\begin{equation*}
\max_{u \in [-1,1]} \Big|\frac{d\paak}{du}\Big|\sim k^2.
\end{equation*}
Hence we have proved that
\begin{equation*}
\lim_{k \to \infty} \left|\paak(t^{(\alpha+1,\beta)}_{k-1,k-1})-\paak\left(\cos\frac{j_{\al+1}}{k}\right)\right| =0,
\end{equation*}
and \eqref{eq:limit} follows.

Since the zeros $t^{(\alpha,\beta)}_{k,k}$ tend to $1$ as $k$ tends to infinity, to prove the theorem it suffices to show that $\lim_{t \to 1} m(t)$ exists. This follows from \eqref{eq:limit} and the following two facts which hold for all $k \geq 2$:
\begin{equation}
\label{eq:fact1}
P^{(\alpha,\alpha)}_k(t^{(\alpha+1,\alpha+1)}_{k-1,k-1}) \leq m(t) \quad \text{for all $t \geq t^{(\alpha+1,\alpha+1)}_{k-1,k-1}$}
\end{equation}
and
\begin{equation}
\label{eq:fact2}
m(t) \leq P^{(\alpha,\alpha)}_{k+1}(t^{(\alpha+1,\alpha)}_{k,k}) \quad \text{for all $t \in [t^{(\alpha+1,\alpha+1)}_{k-1,k-1},t^{(\alpha+1,\alpha+1)}_{k,k}]$.}
\end{equation}
Fact \eqref{eq:fact1} follows from \eqref{ee1} and \cite[(6.4.19)]{AAR}. For establishing fact~\eqref{eq:fact2} we argue as follows: As in the proof of Proposition~\ref{prop:tookmkm}, we use~\eqref{e5} to show that $P^{(\alpha+1,\alpha)}_k$ has exactly one zero in the interval $[t^{(\alpha+1,\alpha+1)}_{k-1,k-1},t^{(\alpha+1,\alpha+1)}_{k,k}]$, namely $t^{(\alpha+1,\alpha)}_{k,k}$. From~\eqref{e6} we then see that $t^{(\alpha+1,\alpha)}_{k,k}$ is the only point in this interval where $P^{(\alpha,\alpha)}_k$ and $P^{(\alpha,\alpha)}_{k+1}$ coincide. Now it follows from the interlacing property that $P^{(\alpha,\alpha)}_k$ is increasing in the interval and that $P^{(\alpha,\alpha)}_{k+1}$ is decreasing in the interval, and we are done.
\end{proof}

\begin{corollary}
We have
\begin{equation*}
\chim(\R^n) \geq 1+ \frac{(j_{\al+1})^\al}{2^\al \Gamma(\al+1) |J_{\al}(j_{\al+1})|},
\end{equation*}
where $\alpha = (n-3)/2$.\qed
\end{corollary}

We use this corollary to derive new lower bounds for $n = 10, \ldots,
24$. We give them in Table~\ref{table:thetable}. For $n = 2, \ldots,
8$ our bounds are worse than the existing ones and for $n = 9$ our
bound is $35$ which is also the best known one.

In fact Oliveira and Vallentin~\cite{OV} show, by different methods,
that the above bound is actually a bound for $\chim(\R^{n-1})$. This then
gives improved bounds starting from $n = 4$. With the use of
additional geometric arguments one can also get a new bound for $n =
3$ in this framework.

\renewcommand{\thetable}{\arabic{section}.\arabic{table}}

\begin{table}[htb]
\begin{tabular}{c|c|c}
    & best lower bound & new lower bound \\
$n$ & previously known for $\chim(\R^n)$ &  for $\chim(\R^n)$ \\
\hline
10 & 45 & 48\\
11 & 56 & 64\\
12 & 70 & 85\\
13 & 84 & 113\\
14 & 102 & 147 \\
15 & 119 & 191\\
16 & 148 & 248 \\
17 & 174 & 319 \\
18 & 194 & 408 \\
19 & 263 & 521 \\
20 & 315 & 662\\
21 & 374 & 839 \\
22 & 526 & 1060\\
23 & 754 & 1336\\
24 & 933 & 1679\\
\end{tabular}
\\[0.3cm]
\caption{Lower bounds for $\chim(\R^n)$.}
\label{table:thetable}
\end{table}

We can also use the corollary to show that our bound is exponential in
the dimension.  To do so we use the inequalities (cf.~\cite[(4.14.1)]{AAR}
and~\cite[Section 15.3, p. 485]{W})
\begin{equation*}
j_{\al+1} > j_{\al} > \al
\end{equation*}
and (cf.~\cite[(4.9.13)]{AAR})
\begin{equation*}
|J_\al(x)|\leq 1/\sqrt{2} 
\end{equation*}
to obtain 
\begin{equation*}
\frac{(j_{\al+1})^\al}{2^\al \Gamma(\al+1) |J_{\al}(j_{\al+1})|} >
\sqrt{2} \frac{\al^\al}{2^\al\Gamma(\al+1)},
\end{equation*}
and with Stirling's formula $\Gamma(\al+1)\sim \al^\al
e^{-\al}\sqrt{2\pi\al}$ we have that the exponential term is
$\left(\frac{e}{2}\right)^\al \sim (1.165)^n$.

\section{Other spaces}
\label{sec:further}

In this section we want to go back to our generalization \eqref{eq:new
theta} of the theta function and discuss its computation in more 
general situations than the one of the graph $G(n,t)$
encountered in Section \ref{sec:generalization2}. We assume that a
compact group $\Gamma$ acts continuously on $V$.  Then, the
computation only depends on the orthogonal decomposition of the space
of $L^2$-functions~\eqref{eq:gen decomp}.

\subsection{Two-point homogeneous spaces}

First, it is worth noticing that all results in
Section~\ref{sec:generalization2} are valid --- one only has to use the
appropriate zonal polynomials and appropriate volumes --- for distance
graphs in infinite, two-point homogeneous, compact metric spaces where
edges are given by exactly one distance. 

If one considers distance graphs in infinite, compact, two-point
homogeneous metric spaces with $s$ distances, then
it is helpful to consider a dual formulation of~\eqref{eq:primal}. It
is an infinite linear programming problem in dimension $s+1$ which in
the case of the unit sphere has the following form:
\begin{equation*}
\begin{split}
\min\big\{ & z_1/\omega_n :  \\
& \quad \text{$z_1 + z_{t_1} + \cdots + z_{t_s} \geq \omega_n^2$},\\
& \quad \text{$z_1 + z_{t_1} \paak(t_1) + \cdots + z_{t_s} \paak(t_s) \geq 0$ for $k = 1, 2, \ldots$}\big\},
\end{split}
\end{equation*}
where $t_1, \ldots, t_s$ are the inner products defining the edges of our
graph.

\subsection{Symmetric spaces}

Next we may consider infinite compact metric spaces $V$ which are not
two-point homogeneous but symmetric. Since the space $L^2(V)$ still
has a multiplicity-free orthogonal decomposition one gets a linear
programming bound, but with the additional complication that one has
to work with multivariate zonal polynomials. The most prominent case
of the Grassmann manifold was considered by the first author in
\cite{Bac} in the context of finding upper bounds for finite codes.

\subsection{General homogeneous spaces}

For the most general case one would have multiplicities $m_k$ in the
decomposition of $L^2(V)$ which is given by the Peter-Weyl Theorem:
\begin{equation}
\label{eq:gen decomp}
L^2(V) = (H_{0,1} \perp \ldots \perp H_{0,m_0}) \perp (H_{1,1} \perp
\ldots \perp H_{1,m_1}) \perp \ldots,
\end{equation}
where $H_{k,l}$ are $\Gamma$-irreducible subspaces which are
equivalent whenever their first index coincides. In this case one uses
Bochner's characterization of the continuous, $\Gamma$-invariant,
positive kernels given in \cite[Section III]{Boc} which yields a true
semidefinite programming problem for the computation of $\vartheta$.

\section{Second generalization}
\label{sec:second}

In this section we first show how our generalization relates to the
theta function of finite subgraphs of $G(n,t)$. We prove that
computing the theta function for any finite subgraph of
$G(n,t)$ does not give a better bound than the one of
Corollary~\ref{cor:lowerboundchi}. For this we introduce a second
generalization of the theta function. Then we show how our second
generalization relates to the linear programming bound of Delsarte.

\subsection{Finite subgraphs}

To compute a bound for the measurable chromatic number of the
graph~$G(n, t)$ we compute~$\vartheta(G(n, t))$, which is an upper
bound for $\alpha(G(n, t))$, and then $\omega_n/\vartheta(G(n, t))$ is
a lower bound for~$\chim(G(n, t))$.

When~$G = (V, E)$ is a finite graph, this approach corresponds to
computing~$\vartheta(G)$ and using $|V| / \vartheta(G)$ as a lower
bound for~$\chi(G)$. However, this is in general not the best bound we
can obtain for~$\chi(G)$ from the theta function. Indeed, for a finite
graph~$G$, the so-called \defi{sandwich theorem} says that
\begin{equation*}
\alpha(G) \leq \vartheta(G) \leq \chi(\overline{G})
\end{equation*}
(Theorem~\ref{th:original} only gives the first inequality, Lov\'asz
\cite[Proof of Corollary 3]{Lov1} gives the second),
where~$\overline{G}$ is the \defi{complement}
of~$G$, the graph with the same vertex set as~$G$ and in which two
vertices are adjacent if and only if they are nonadjacent in~$G$.

Moreover, for a finite graph $G = (V, E)$, we have
\begin{equation}
\label{eq:thetaprod}
\vartheta(G) \vartheta(\overline{G}) \geq |V|
\end{equation}
(cf.~Lov\'asz~\cite[Corollary~2]{Lov1}). For some graphs (e.g., stars),
this inequality is strict, hence in these
cases~$\vartheta(\overline{G})$ would provide us with a better lower
bound for~$\chi(G)$ than $|V| / \vartheta(G)$ would. But when~$V$ is
homogeneous we actually have equality in~\eqref{eq:thetaprod}
(cf.~Lov\'asz~\cite[Theorem 8]{Lov1}). In this case, both bounds
for~$\chi(G)$ coincide.

Something similar happens for our infinite distance graph~$G(n,
t)$. The complement of~$G(n, t)$ is the graph in which any two
distinct points on the unit sphere whose inner product is not~$t$ are
adjacent. We cannot use our generalization of the theta function to
define~$\vartheta(\overline{G(n, t)})$. However, we may use a
different (and for finite graphs, equivalent) definition
of~$\vartheta$ (cf.~Lov\'asz~\cite[Theorem~3]{Lov1}), which for a
finite graph~$G = (V, E)$ is 
\begin{equation}
\label{eq:theta2}
\begin{split}
  \vartheta(\overline{G}) = \min\big\{ \lambda : \quad & \text{$K \in \R^{V
\times V}$ is positive semidefinite},\\ 
& \text{$K(x,x) = \lambda - 1$ for all $x \in V$,}\\ 
& \text{$K(x,y) = -1$ if $\{x,y\} \in E$}\big\}.
\end{split}
\end{equation}

The generalization of this definition, applied to~$\overline{G(n, t)}$
and with the symmetry taken into account, is described below. We
choose to write~$\overline{\vartheta}(G(n, t))$ instead
of~$\vartheta(\overline{G(n, t)})$ to emphasize that the two ways to
define the theta function are \textsl{not} equivalent for our infinite
graph. So we have
\begin{equation}
\label{eq:theta2gen}
\begin{split}
  \overline{\vartheta}(G(n,t)) = \min\big\{ \lambda : \quad & \text{$K \in \Ccal(S^{n-1} \times S^{n-1})$ is positive},\\ 
  & \text{$K$ is invariant under $\ort(\R^n)$},\\ 
& \text{$K(x,x) = \lambda - 1$ for all $x \in S^{n-1}$,}\\ 
& \text{$K(x,y) = -1$ if $x \cdot y = t$}\big\}.
\end{split}
\end{equation}

By decomposing the kernel~$K$ with the help of the Jacobi polynomials
as done in Section~\ref{sec:generalization2}, we may compute the
optimal value of the optimization problem~\eqref{eq:theta2gen}, and in
doing so we find out that
\begin{equation*}
\vartheta(G(n, t)) \overline{\vartheta}(G(n, t)) = \omega_n,
\end{equation*}
so that we have the analogue of $\vartheta(G) \vartheta(\overline{G})
= |V|$ for our infinite distance graph on the unit sphere.

This also provides us with the connection to the theta function of
finite subgraphs of~$G(n, t)$ claimed in
Section~\ref{sec:generalization2}. If~$H = (V, E)$ is a finite
subgraph of~$G(n, t)$, then~$\vartheta(\overline{H})$ provides a lower
bound for~$\chi(H)$, which in turn is a lower bound for~$\chim(G(n,
t))$. It could be that for some finite subgraph~$H$ of~$G(n, t)$ this
lower bound would be better than the one provided by~$\vartheta(G(n,
t))$.  This is, however, not the case. Indeed, if~$K$ is an optimal
solution for~\eqref{eq:theta2gen}, the restriction of~$K$ to~$V \times
V$ is a feasible solution to the optimization
problem~\eqref{eq:theta2} defining~$\vartheta(\overline{H})$, hence
$\vartheta(\overline{H}) \leq \overline{\vartheta}(G(n, t))$, which is
our bound for~$\chim(G(n, t))$.

\subsection{Delsarte's linear programming bound}

The second generalization $\overline{\vartheta}$ of the theta function
is closely related to the linear programming bound for finite codes
established by Delsarte in \cite{Del} and put into the framework of
group representations, which we use here, by Kabatiansky and
Levenshtein in \cite{KL}. Here we devise an explicit connection
between these two bounds. The connection between the linear
programming bound and the theta function was already observed by
McEliece, Rodemich,  Rumsey Jr.\ in \cite{MRR} and
independently by Schrijver in \cite{S} in the case of finite
graphs.

Consider the graph on the unit sphere where two distinct points are
adjacent whenever their inner product lies in the open interval
$[-1,t]$. We denote this graph by $G(n,[-1,t])$. Stable sets in the
\textsl{complement} of this graph are finite and consist of points on the unit
sphere with minimal angular distance $\arccos t$.

Now the second generalization \eqref{eq:theta2} applied to $G(n,[-1,t])$ is
\begin{equation}
\label{eq:theta2gen gnt1}
\begin{split}
  \overline{\vartheta}(G(n,[-1,t])) = \inf\big\{ \lambda : \quad & \text{$K \in \Ccal(S^{n-1} \times S^{n-1})$ is positive},\\ 
  & \text{$K$ is invariant under $\ort(\R^n)$},\\ 
& \text{$K(x,x) = \lambda - 1$ for all $x \in S^{n-1}$,}\\ 
& \text{$K(x,y) = -1$ if $x \cdot y \in [-1,t]$}\big\}.
\end{split}
\end{equation}
We safely write $\inf$ instead of $\min$ here because we do not know
if the infimum is attained. 

\begin{proposition}
\label{prop:theta2bound}
Let $C \subseteq S^{n-1}$ be a subset of the unit sphere such that every pair
of distinct points in $C$ has inner product lying in $[-1,t]$. Then its cardinality is at
most $\overline{\vartheta}(G(n,[-1,t]))$.
\end{proposition}

\begin{proof}
Let $K$ be a kernel satisfying the conditions in
\eqref{eq:theta2gen gnt1}. Then, by the positivity of the continuous kernel
$K$ it follows that
\begin{equation*}
0 \leq \sum_{(c,c') \in C^2} K(c,c')
= \sum_{c} K(c,c) + \sum_{c \neq c'} K(c,c')
\leq |C| K(c,c) - |C|(|C| - 1),
\end{equation*}
so that $|C| - 1 \leq K(c,c)$ and we are done.
\end{proof}

We finish by showing how the original formulation of the linear
programming bound can be obtained from \eqref{eq:theta2gen gnt1}. Using
Schoenberg's characterization~\eqref{eq:schoenberg} the semidefinite
programming problem~\eqref{eq:theta2gen gnt1} simplifies to the linear
programming problem
\begin{equation*}
\begin{split}
\inf\big\{ \lambda : \quad & \text{$f_0 \geq 0, f_1 \geq 0, \ldots$},\\ 
& \text{$\sum_{k=0}^\infty f_k \paak(1) = \lambda-1$},\\
& \text{$\sum_{k=0}^\infty f_k \paak(u) = -1$ for all $u \in [-1,t]$}
\big\}.
\end{split}
\end{equation*}
We can strengthen it by requiring $\sum_{k=0}^\infty f_k \paak(u) \leq
-1$ for all $u \in [-1,t]$. By restricting $f_0 = 0$ the infimum is
not effected.  Then, after simplification, we get the linear
programming bound (cf. \cite{DGS}, \cite{KL}).
\begin{equation*}
\begin{split}
\inf\{1 + \sum_{k = 1}^\infty f_k : \quad & \text{$f_1 \geq 0, f_2 \geq 0, \ldots$},\\ 
& \text{$\sum_{k=1}^\infty f_k \paak(u) \leq -1$ for all $u \in [-1,t]$}
\big\}.
\end{split}
\end{equation*}
By Proposition~\ref{prop:theta2bound} it gives an upper bound for the
maximal number of points on the unit sphere with minimal angular
distance $\arccos t$.

\section*{Acknowledgements}

We thank Dion Gijswijt, Gil Kalai, Tom Koornwinder, Pablo Parrilo, and Lex Schrijver
for their helpful comments.

\end{document}